# An extension of the Artin-Mazur theorem

By Vadim Yu. Kaloshin

*Dedicated to the memory of my grandfather Meyer Levich*

## 1. Introduction

Let $C^r(M, M)$ be the space of $C^r$ mappings of a compact manifold $M$ into itself with the uniform $C^r$-topology and $\text{Diff}^r(M)$ be the space of $C^r$ diffeomorphisms of $M$ with the same topology. It is well-known that $\text{Diff}^r(M)$ is an open subset of $C^r(M, M)$. For a map $f \in C^r(M)$, consider the number of *isolated* periodic points of period $n$ (i.e. the number of isolated fixed points of $f^n$)

(1) $$P_n(f) = \#\{ \text{ isolated } \ x \in M : x = f^n(x)\}.$$

In 1965 Artin and Mazur [AM] proved the following:

THEOREM 1. *There exists a dense set $\mathcal{D}$ in $C^r(M, M)$ such that for any map $f \in \mathcal{D}$ the number $P_n(f)$ grows at most exponentially with $n$, i.e. for some number $C > 0$*

(2) $$P_n(f) \leq \exp(Cn) \quad \text{forall} \ \ n \in \mathbb{Z}_+.$$

Notice that this theorem does not exclude the possibility that a mapping $f$ in $\mathcal{D}$ has a curve $\gamma$ of periodic points; i.e., for all $x \in \gamma$, $f^n(x) = x$ for some $n \in \mathbb{Z}_+$, because in this case $\gamma$ consists of nonisolated periodic points of period $n$ (see the last part of Theorem 6 for this nonisolated case).

*Definition* 2. We call a mapping (resp. diffeomorphism) $f \in C^r(M, M)$ (resp. $f \in \text{Diff}^r(M)$) an Artin-Mazur mapping (resp. diffeomorphism) or simply an A-M mapping (resp. diffeomorphism) if $P_n(f)$ grows at most exponentially fast.

Artin and Mazur [AM] posed the following problem: *What can be said about the set of* A-M *mappings with only transversal periodic orbits*? Recall that a periodic orbit of period $n$ is called *transversal* if the linearization $df^n$ at this point has for an eigenvalue no $n^{\text{th}}$ roots of unity. Notice that a hyperbolic periodic point is always transversal, but not vice versa.

In what follows we consider not the whole space $C^r(M, M)$ of mappings of $M$ into itself, but only its open subset $\text{Diff}^r(M)$. The main result of this paper is the following.

MAIN THEOREM. *Let $1 \leq r < \infty$. Then the set of A-M diffeomorphisms with only hyperbolic periodic orbits is dense in the space $\text{Diff}^r(M)$.*

This theorem says that A-M diffeomorphisms which satisfy part of the Kupka-Smale condition form a dense set in $\text{Diff}^r(M)$. Recall that a diffeomorphism is called a *Kupka-Smale* (or K-S) *diffeomorphism* if all its periodic points are hyperbolic and all associated stable and unstable manifolds intersect one another transversally. The Kupka-Smale theorem says that K-S diffeomorphisms form a residual set (see e.g. [PM]). The natural question is whether the intersection of A-M and K-S diffeomorphisms is dense. The answer is not easy, because methods of the proof of both theorems are of completely different nature and cannot be applied simultaneously. If one omits the condition on transversality of stable and unstable manifolds, then the Main Theorem says that the intersection of A-M and K-S diffeomorphisms is dense.

A residual set in a finite-dimensional space can have measure zero. Therefore, the Kupka-Smale theorem does not imply that "almost every" diffeomorphism is a K-S diffeomorphism. In loose terms, a set $P \subset \text{Diff}^r(M)$ is called *prevalent* if for a generic finite parameter family $\{f_\epsilon\}_{\epsilon \in \text{Ball}}$, the property $f_\epsilon \in P$ holds for almost every parameter value. For a discussion of prevalence in linear space see also [HSY]. Finally, in [K2] it is proven that K-S diffeomorphisms form a prevalent set which is a stronger (in some sense) statement than the classical one. In Section 3 we present some additional results about the set of A-M diffeomorphisms. It turns out that the set of A-M diffeomorphisms is *not* residual.

## 2. A Proof of the Main Theorem

Consider a $C^r$ diffeomorphism $f : M \to M$. We shall approximate $f$ by an A-M diffeomorphism with only hyperbolic periodic orbits. The proof of it consists of two steps.

*Step* 1. Reduction to a problem for polynomial maps.    Using the Whitney embedding theorem, embed $M$ into $\mathbb{R}^N$ for $N = 2 \dim M + 1$. Denote by $T$ a tube neighborhood of $M$. For any fixed $r \in \mathbb{Z}_+$ one can extend $f : M \to M$ to a diffeomorphism $F : T \to T$ of the tube neighborhood $T$ strictly into inself such that $F$ restricted to $M$ coincides with $f$. If $F$ contracts along the directions transverse to $M$ sufficiently strongly, then by the Sacker theorem [Sa] (see also [Fe], [HPS]) each diffeomorphism $\tilde F : T \to T$ which is $C^r$-close to $F$ has a $C^r$ smooth invariant manifold $\tilde M$ which is $C^r$-close to $M$. Denote by $\pi : \tilde M \to M$ a diffeomorphism from $\tilde M$ to $M$ which can be obtained by projection along the normal to $M$ directions. Then $\tilde f = \pi \circ \tilde F|_{\tilde M} \circ \pi^{-1} : M \to M$ is a diffeomorphism which is $C^r$-close to $f$. By the Weierstrass approximation theorem one can approximate a diffeomorphism $F : T \to T$ of an open set $T$ in the Euclidean space $\mathbb{R}^N$ into itself by a polynomial diffeomorphism $\tilde F = P|_T : T \to T$. Notice that if $\tilde F$ has only hyperbolic periodic orbits, then the induced diffeomorphism $\tilde f = \pi \circ \tilde F|_{\tilde M} \circ \pi^{-1} : M \to M$ also has only hyperbolic periodic orbits.

We shall prove that, indeed, one can approximate any diffeomorphism $F : T \to T$ by a polynomial diffeomorphism $\tilde F = P|_T : T \to T$ which has only hyperbolic periodic orbits.

Let $D \in \mathbb{Z}_+$. Denote by $A_N^D$ the space of (vector-)polynomials $P : \mathbb{R}^N \to \mathbb{R}^N$ of degree at most $D$. If $\mu = \mu(N, D) = \#\{\alpha \in \mathbb{Z}_+^N : |\alpha| \leq D\}$, then $A_N^D$ is isomorphic to $\mathbb{R}^\mu$. Consider $A_N^D$ with the Lebesgue measure on it.

*Step 2.* For any $D \in \mathbb{Z}_+$, almost every polynomial $P : \mathbb{R}^N \to \mathbb{R}^N$ from $A_N^D$ has only hyperbolic periodic orbits and their number grows at most exponentially.

The second part of this statement is easy provided that the first is true. Indeed, fix $k \in \mathbb{Z}_+$, $k > 0$ and consider the system

$$P(x_1) - x_2 = 0, \ P(x_2) - x_3 = 0, \ldots, \ P(x_k) - x_1 = 0.$$

This system has $Nk$ equations, each of them of degree at most $D$. By the Bezout theorem the number of isolated solutions is at most $D^{kN}$. If all periodic points are hyperbolic, then they are all isolated and this completes the proof.

Fix $k \in \mathbb{Z}_+$, $k > 0$. Let $\alpha = (\alpha_1, \ldots, \alpha_N) \in \mathbb{Z}_+^N$ be a multiindex, $|\alpha| = \sum_i \alpha_i$. Fix a coordinate system in $\mathbb{R}^N$ so one can write each polynomial $P(a, \cdot) : \mathbb{R}^N \to \mathbb{R}^N$ from $A_N^D$ in the form

(3)
$$P(a, x) = \sum_{|\alpha| \leq D} a_\alpha x^\alpha, \ \text{where} \ a = (\{a_\alpha\}_{|\alpha| \leq D}) \in \mathbb{R}^\mu,$$

$$x = (x_1, \ldots, x_N) \in \mathbb{R}^N, \ \text{and} \ x^\alpha = x_1^{\alpha_1} \ldots x_N^{\alpha_N}.$$

LEMMA 1. *Let $\lambda_0 \in \mathbb{C}$ and $|\lambda_0| = 1$. For any $D \in \mathbb{Z}_+$, almost every polynomial $P : \mathbb{R}^N \to \mathbb{R}^N$ from $A_N^D$ has no periodic orbits with the eigenvalue $\lambda_0$.*

Denote the $k$-fold composition $P(a, \cdot) \circ \ldots \circ P(a, \cdot) : \mathbb{R}^N \to \mathbb{R}^N$ by $P^{(k)}(a, \cdot)$, the linearization matrix of the map $P^{(k)}(a, \cdot)$ at a point $x$ by $d_x(P^{(k)})(a, x)$, and the $N \times N$ identity matrix by $\mathrm{Id}_N$. Let $\lambda \in \mathbb{C}$ be a complex number.

Denote $D(a, \lambda, x) = \det (d_x (P^{(k)}) (a, x) - \lambda \, \mathrm{Id}_N)$. Every periodic orbit of period $k$, which has an eigenvalue $\lambda$, satisfies the following system:

(4)
$$\begin{cases} P^{(k)}(a, x) - x = 0, & x = (x_1, \ldots x_N) \in \mathbb{R}^N \\ D(a, \lambda, x) = 0, & a \in \mathbb{R}^\mu. \end{cases}$$

The general goal is to prove that for a "generic" choice of coefficients $a \in \mathbb{R}^\mu$ of $P(a, \cdot)$ this system has no solutions satisfying the condition $|\lambda| = 1$ or there is no nonhyperbolic periodic orbit of period $k$. First, we prove that a "generic" choice of coefficients $a \in \mathbb{R}^\mu$ of $P(a, \cdot)$ has no periodic points with the eigenvalue $\lambda = \lambda_0$.

Notice that the system (4) including the condition $\lambda = \lambda_0$ (or $|\lambda| = 1$) consists of $N + 2$ equations and for each value $a$ only $N + 1$ variables $x_1, \ldots, x_N, \lambda$. It might be clear intuitively that for a "generic" $a \in \mathbb{R}^\mu$ there is no solution, because the number of equations is more than the number of variables. To prove it rigorously for $\lambda = \lambda_0$ (or $|\lambda| = 1$) we shall apply elimination theory.

2.1. *Elimination theory.* Let $\mathbb{C}^m$ denote the $m$-dimensional complex space $z = (z_1, \ldots, z_m) \in \mathbb{C}^m$, $m \in \mathbb{Z}_+$. A set $V$ in $\mathbb{C}^m$ is called *a closed algebraic set* in $\mathbb{C}^m$ if there is a finite set of polynomials $F_1, \ldots F_s$ in $z_1, \ldots, z_m$ such that

$$V(F_1, \ldots, F_s) = \{(z_1, \ldots, z_m) \in \mathbb{C}^m | \ F_j(z_1, \ldots, z_m) = 0, \ 1 \leq j \leq s\}.$$

One can define a topology in $\mathbb{C}^m$, called the *Zariski topology*, whose closed sets are closed algebraic sets in $\mathbb{C}^m$. This, indeed, defines a topology, because the set of closed algebraic sets is closed under a finite union and an arbitrary intersection. Sometimes, closed algebraic sets are also called Zariski closed sets.

*Definition* 3. A subset $S$ of $\mathbb{C}^m$ is called constructible if it is in the Boolean algebra generated by the closed algebraic sets; or equivalently if $S$ is a disjoint union $T_1 \cup \ldots \cup T_k$, where $T_i$ is locally closed, i.e. $T_i = T'_i - T''_i$, where $T'_i$ is a closed algebraic set and $T''_i \subset T'_i$ is a smaller closed algebraic set.

One of the main results of elimination theory is the following:

THEOREM 4 ([Mu, Ch. 2.2]). *Let $V \subset \mathbb{C}^\mu \times \mathbb{C}^N$ be a constructible set and $\pi : \mathbb{C}^\mu \times \mathbb{C}^N \to \mathbb{C}^\mu$ be the natural projection. Then $\pi(V) \subset \mathbb{C}^\mu$ is a constructible set.*

*Remark* 1. An elementary description of elimination theory can be found in books by Jacobson [J] and van der Waerden [W].

2.2. *Proof of Lemma* 1 *or application of elimination theory to the system* (4). Put $\lambda = \lambda_0$ and consider the system (4) for $(a; x) \in \mathbb{C}^\mu \times \mathbb{C}^N$. Then it defines a closed algebraic set $V_k(\lambda_0) \subset \mathbb{C}^\mu \times \mathbb{C}^N$. By Theorem 4 the natural projection $\pi : \mathbb{C}^\mu \times \mathbb{C}^N \to \mathbb{C}^\mu$ of $V_k(\lambda_0)$, namely $\pi(V_k(\lambda_0)) \subset \mathbb{C}^\mu$, is a constructible set. The only thing left to show is that $\pi(V_k(\lambda_0)) \neq \mathbb{C}^\mu$ and has a positive codimension. Recall that $|\lambda_0| = 1$ and $\mathbb{C}^\mu$ is the space of complex coefficients of a polynomial of degree $D$.

LEMMA 2. *Let $\mathbb{R}^\mu$ be naturally embedded into $\mathbb{C}^\mu$. Then there is an open set $U \subset \mathbb{C}^\mu$ such that $U \cap \mathbb{R}^\mu \neq \emptyset$ and for any $a \in U$ the corresponding polynomial $P(a, \cdot) : \mathbb{C}^N \to \mathbb{C}^N$ of degree $0$ has exactly $D^{Nk}$ periodic points of period $k$ and all of them are hyperbolic. Therefore, $U \cap \pi(V_k(\lambda_0)) = \emptyset$ and $\pi(V_k(\lambda_0))$ is a constructible set of a positive codimension in $\mathbb{C}^\mu$.*

In the case $D = 1$ the statement of the proposition is obvious.

*Proof.* Consider the homogeneous polynomial $P(a^*, \cdot) : \mathbb{C}^N \to \mathbb{C}^N$ of degree $D$ given by $P(a^*, \cdot) : (z_1, \ldots, z_N) \mapsto (z_1^D \ldots, z_N^D)$. It is easy to see that $P$ has exactly $D^{Nk}$ periodic points of period $k$ all of which are hyperbolic. Hyperbolicity of periodic points of period $k$ of $P$ implies that any polynomial mapping $\tilde{P}$, which is a small perturbation of $P$, has at least $D^{Nk}$ hyperbolic points of period $k$ close to those of $P$. Since, a polynomial of degree $D$ has at most $D^{Nk}$ periodic point of period $k$, Bezout's theorem implies that there are exactly $D^{Nk}$ periodic point of period $k$, and the set of periodic points of period $k$ has no components of positive dimension (see e.g. [Sh, Ch. 4.2]). Thus, there is a neighborhood $U \subset \mathbb{C}^\mu$ of $a^*$ such that for any $a \in U$ the corresponding polynomial $P(a, \cdot) : \mathbb{C}^N \to \mathbb{C}^N$ has only hyperbolic periodic points of period $k$ and the definition of $V_k(\lambda_0)$ implies that $U \cap \pi(V_k(\lambda_0)) = \emptyset$. Since $\pi(V_k(\lambda_0))$ is constructible and does not intersect an open set, this implies that $\pi(V_k(\lambda_0))$ has dimension less than $\mu$ in $\mathbb{C}^\mu$. This completes the proof of Lemma 2.

By Proposition 2 the restriction $\pi(V_k(\lambda_0)) \cap \mathbb{R}^\mu$ has positive codimension and, therefore, measure zero in $\mathbb{R}^\mu$. Thus, almost every polynomial $P(a, \cdot)$ from $A_N^D = \mathbb{R}^\mu$ has no periodic points of period $k$ with the eigenvalue $\lambda_0$. Intersection over all $k \in \mathbb{Z}_+$ gives that the same is true for all periods. This completes the proof of Lemma 1.

2.3. *Completion of the proof of Step* 2 *of the Main Theorem.* Consider the system (4) for $(a, \lambda; x) \in \mathbb{C}^\mu \times \mathbb{C} \times \mathbb{C}^N$. It defines a closed algebraic set, denoted by $V_k \subset \mathbb{C}^\mu \times \mathbb{C} \times \mathbb{C}^N$. By Theorem 4 the natural projection $\pi : \mathbb{C}^\mu \times \mathbb{C} \times \mathbb{C}^N \to \mathbb{C}^\mu \times \mathbb{C}$ of $V_k$, namely $S_k = \pi(V_k) \subset \mathbb{C}^\mu \times \mathbb{C}$, is a constructible set.

Consider natural projections $\pi_1 : \mathbb{C}^\mu \times \mathbb{C} \to \mathbb{C}^\mu$ and $\pi_2 : \mathbb{C}^\mu \times \mathbb{C} \to \mathbb{C}$. It follows from Proposition 2 that $S_k$ has dimension $\mu$. Indeed, with the notation of Proposition 2, the projection $\pi_1(S_k) = W_k$ contains an open set $U \subset \mathbb{C}^\mu$ and $S_k$ does not intersect a neighborhood of $U \times \{\lambda : |\lambda| = 1\} \subset \mathbb{C}^\mu \times \mathbb{C}$.

By Theorem 4, $\pi_1(S_k) = W_k$ is constructible and by Proposition 2, $\dim W_k = \mu$. By Sard's lemma for algebraic sets ([Mu, Ch. 3.A]) there exists a proper algebraic set $\Sigma_k \subset S_k$ such that the map $\pi_1$ is restricted to $\tilde{S}_k = S_k \setminus \Sigma$; namely, $p_k = \pi_1|_{\tilde{S}_k} : \tilde{S}_k \to \mathbb{C}^\mu$ has no critical points. Thus, outside of some smaller closed algebraic set $\Sigma'_k \subset \mathbb{C}^\mu$ the map $p_k : \tilde{S}_k \to \mathbb{C}^\mu \setminus \Sigma'_k$ is locally invertible.

Recall that our goal is to show that $Z = \pi_1(S_k \cap \{\lambda : |\lambda| = 1\}) \cap \mathbb{R}^\mu$ has measure zero in $A_N^D = \mathbb{R}^\mu$. It is sufficient to prove this locally.

Let $a \in \mathbb{R}^\mu \setminus \Sigma'_k$ and $U \subset \mathbb{C}^\mu \setminus \Sigma'_k$ be a neighborhood of $a$. By construction the map $p_k : \tilde{S}_k \to \mathbb{C}^\mu \setminus \Sigma'_k$ is locally invertible, so the preimage $p_k^{-1}(U)$ consists of a finite disjoint union of open sets $\cup_{j \in J} U_j \subset S_k$. Thus, one can define a finite collection of analytic functions $\{\lambda_{k,j} = \pi_2 \circ p_{k,j}^{-1} : U \to \mathbb{C}\}_{j \in J}$, where $p_{k,j}^{-1} : U \to U_j$ is the inverse of the restriction $p_k|_{U_j} : U_j \to U$. We need to show that

$$\cup_{j \in J} \left( \lambda_{k,j}^{-1}(\{\lambda : |\lambda| = 1\}) \right) \cap \mathbb{R}^\mu$$

has measure zero. If for some $j \in J$ the function $\lambda_{k,j} : U \to \mathbb{C}$ is equal to a constant $\lambda$, then Proposition 2 implies $|\lambda| \neq 1$ and the preimage

$$\lambda_{k,j}^{-1}(\{\lambda : |\lambda| = 1\}) = \emptyset.$$

If for some $j \in J$ the function $\lambda_{k,j} : U \to \mathbb{C}$ is not constant, then the set $\lambda_{k,j}^{-1}(\{\lambda : |\lambda| = 1\}) \cap \mathbb{R}^\mu$ is a real analytic set of a positive codimension and, therefore, has measure zero. It follows e.g. from the fact that a real analytic set can be stratified (see e.g. [H] or [GM]), i.e., in particular, can be decomposed into at most a countable union of semianalytic manifolds. Each semianalytic manifold must have a positive codimension and, therefore, measure zero. This implies that for almost every $a \in \mathbb{R}^\mu$ the system (4) has no solutions for $|\lambda| = 1$. Since, $k$ is arbitrary, this completes the proof of step 2.

Let us complete the proof of the Main Theorem. Application of Step 1 shows that a diffeomorphism $f : M \to M$ can be extended to a tube neighborhood $T$ of $M$, $F : T \to T$, and that it is sufficient to approximate $F$ by a diffeomorphism $\tilde{F} : T \to T$ which has only hyperbolic periodic points. By the Weierstrass approximation theorem, $F$ can be approximated by a polynomial diffeomorphism $\tilde{F} = P|_T : T \to T$. Since, in the space of polynomial maps, of any degree $D$, polynomial maps with only hyperbolic periodic points form a full measure set, one can choose $\tilde{F} = P|_T : T \to T$ which has only hyperbolic periodic points. If $\tilde{F} : T \to T$ has only hyperbolic periodic points, then its restriction to an invariant manifold also has only hyperbolic periodic points. This completes the proof of the Main Theorem. □

*Remark* 2. In order to give a positive answer to the Artin-Mazur question stated in the introduction it is sufficient to use only Step 1 and Lemma 1 of the above proof.

## 3. Generic diffeomorphisms with superexponential growth of the number of periodic points

Let us formulate two important questions related to the growth of the number of periodic points.

In [AM] Artin-Mazur introduced the *dynamical $\zeta_f$-function* defined by $\zeta_f(z) = \exp\left(\sum_{n=1}^{\infty} P_n(f)\frac{z^n}{n}\right)$. Recall that $P_n(f)$ denotes the number of isolated periodic orbits of $f$ of period $n$. In general, the dynamical $\zeta_f$-function is a formal power series in $z$, which may not have a positive radius of convergence. The A-M diffeomorphisms are characterized by the property that the radius of convergence is positive. It is well-known that the dynamical $\zeta_f$-function of a diffeomorphism $f$ satisfying Axiom A has an analytic continuation to a rational function (see [Ba]).

In 1967 Smale [Sm] posed the following question (Problem 4.5, p. 765):

*Is a dynamical $\zeta_f$-function generically rational (i.e. is $\zeta_f$ rational for a residual set of $f \in \mathrm{Diff}^r(M)$)?*

In [Si] it is shown that for the 3-dimensional torus there is no residual set on which the dynamical $\zeta_f$-function is rational. It turns out that for manifolds of dimension greater than or equal to 2 and $C^r$ diffeomorphisms with $r \geq 2$ it is not even analytic in any disk around $z = 0$ (see below). Recall that a subset of a topological space is called residual if it contains a countable intersection of open dense subsets. The usual terminology is to say "a property of points in a Baire space is (topologically) generic" if the set of points which satisfies this property is residual.

Finally, in 1978 R. Bowen asked the following question in [Bo]:

Let $h(f)$ denote the topological entropy of $f$. *Is the property that*

$$h(f) = \limsup_{n \to \infty} \log P_n(f)/n$$

*generic with respect to the $C^r$ topology?*

It turns out that the two questions above can be answered simultaneously for $C^r$ diffeomorphisms with $2 \leq r < \infty$. Below we formulate a theorem of Gonchenko-Shil'nikov-Turaev [GST] which implies the following:

THEOREM 5 ([K1]). *The A-M property is not $C^r$-generic.*

COROLLARY 1. *The property of having a convergent $\zeta_f(z)$ function is not $C^r$-generic, nor is the equation $h(f) = \limsup_{n \to \infty} \log P_n(f)/n$.*

The first part follows from the fact that the A-M property is equivalent to $\zeta_f(z)$ having a positive radius of convergence. To prove the second part notice that the topological entropy for any $C^r$ ($r \geq 1$) diffeomorphism $f$ of a compact manifold is always finite (see e.g. [HK]), whereas if $f$ does not have the A-M property, then $\limsup_{n\to\infty} \log P_n(f)/n = \infty$.

Since an Axiom A diffeomorphism is an A-M diffeomorphism, we need to analyze the complement to the set of Axiom A diffeomorphisms in the space $\text{Diff}^r(M)$. Notice that an example of a diffeomorphism with an arbitrarily fast growing number of periodic orbits is given in [RG]. Now we describe a "bad" domain, where the A-M property fails to be topologically generic.

Let us describe a class of open sets in $\text{Diff}^r(M)$, where the A-M property fails to be topologically generic.

In 1970 Newhouse found an open set in the space of $C^r$ diffeomorphisms $\text{Diff}^r(M)$, where diffeomorphisms exhibiting homoclinic tangencies, defined below, are dense [N]. Such a domain is called *a Newhouse domain* $\mathcal{N} \subset \text{Diff}^r(M)$.

THEOREM 6. *Let* $2 \leq r < \infty$. *Let* $\mathcal{N} \subset \text{Diff}^r(M)$ *be a Newhouse domain. Then for an arbitrary sequence of positive integers* $\{a_n\}_{n=1}^\infty$ *there exists a residual set* $\mathcal{R}_a \subset \mathcal{N}$, *depending on the sequence* $\{a_n\}_{n=1}^\infty$, *with the property that* $f \in \mathcal{R}_a$ *implies that*

$$\limsup_{n\to\infty} P_n(f)/a_n = \infty.$$

*Moreover, there is a dense set* $\mathcal{D}$ *in* $\mathcal{N}$ *such that any diffeomorphism* $f \in F$ *has a curve of periodic points.*

This is another theorem which follows from the theorem of Gonchenko-Shil'nikov-Turaev, which will be discussed in the next section.

In such a domain, Newhouse exhibited a residual set of diffeomorphisms with infinitely many distinct sinks [N], [R], and [PT]. Now it is known as *Newhouse's phenomenon*. Theorem 6 is similar to Newhouse's phenomenon in the sense that a strange phenomenon occurs on a residual set.

Theorem 5 is a corollary of the first part of Theorem 6. To see this fix the sequence $a_n = n^n$ and denote by $\mathcal{R}_a$ a set from Theorem 6 corresponding to this sequence. Assume that A-M diffeomorphisms form a residual set; then this set must intersect with $\mathcal{R}_a$, which is a contradiction.

It seems that based on Newhouse's phenomenon in the space $\text{Diff}^1(M)$ with the $C^1$-topology, where $\dim M \geq 3$, found by Bonnati and Diaz [BD], one can extend Theorems 5 and 6 to the case $r=1$ and $\dim M \geq 3$. The problem with the straightforward generalization is that the proof of the Gonchenko-Shilnikov-Turaev (GST) theorem is two-dimensional in an essential way. To generalize the GST theorem to the three-dimensional case one needs either to find an invariant two-dimensional surface and use the two-dimensional proof or find another proof.

Analogs of Theorems 5 and 6 can be formulated for the case of vector fields on a compact manifold of dimension at least 3. Reduction from the case of diffeomorphisms to the case of vector fields can be done using the standard suspension of a vector field over a diffeomorphism [PM].

Newhouse showed that a Newhouse domain exists under the following hypothesis:

Let a diffeomorphism $f \in \mathrm{Diff}^r(M)$ have a saddle periodic orbit $p$. Suppose stable $W^s(p)$ and unstable $W^u(p)$ manifolds of $p$ have a quadratic tangency. Such a diffeomorphism $f$ is called *a diffeomorphism exhibiting a homoclinic tangency*. Then, arbitrarily $C^r$-close to $f$ in $\mathrm{Diff}^r(M)$, there exists *a Newhouse domain*. In particular, it means that by a small $C^r$-perturbation of a diffeomorphism $f$ with a homoclinic tangency one can generate arbitrarily quick growth of the number of periodic orbits.

On this account we would also like to mention the following conjecture, which is due to Palis [PT], about the space of diffeomorphisms of 2-dimensional manifolds:

CONJECTURE. *If* $\dim M = 2$, *then every diffeomorphism* $f \in \mathrm{Diff}^r(M)$ *can be approximated by a diffeomorphism which is either hyperbolic or exhibits a homoclinic tangency.*

This conjecture is proved for approximation in the $C^1$ topology [PS]. If this conjecture is true in the $C^r$- topology, then in the complement to the set of hyperbolic diffeomorphisms a generic diffeomorphism has arbitrarily quick growth of the number of periodic orbits.

Unfolding of homoclinic tangencies is far from being understood. In [GST] the authors describe the following important result: *there does not exist a finite number of parameters to describe all bifurcations occurring next to a homoclinic tangency* (see §2, Cor. 2 for details). This implies that the complete description of bifurcations of diffeomorphisms with a homoclinic tangency is impossible.

Now we consider diffeomorphisms of a two-dimensional compact manifold $M$. We prove Theorem 6 for two-dimensional manifolds first and then, using the Sacker theorem, we deduce a proof of Theorem 6 for arbitrary $\dim M > 2$ from the two-dimensional case.

## 4. Degenerate periodic orbits in a Newhouse domain and the Gonchenko-Shilnikov-Turaev Theorem [GST]

Assume that a $C^r$ diffeomorphism $f$ exhibits a homoclinic tangency. By the Newhouse theorem [N], in each $C^r$-neighborhood of a diffeomorphism $f$ exhibiting a homoclinic tangency there exists a Newhouse domain.

Let us define a degenerate periodic point of order $k$ or a $k$-degenerate periodic point. Sometimes, it is also called a *saddlenode periodic orbit of multiplicity $k + 1$*.

*Definition* 7. Let $f$ be a $C^s$ diffeomorphism of a 2-dimensional manifold having a periodic orbit $p$ of period $m$. A periodic point $p$ is called *k-degenerate*, where $k < s$, if the linear part of $f^m$ at point $p$ has a multiplier $\nu = 1$ while the other multiplier is different in absolute value from the unit. A *restriction* of $f$ to the central manifold in some coordinate system can be written in the form

(5) $$x \mapsto x + l_{k+1}x^{k+1} + o(x^{k+1}).$$

Let $s > r$. Then $C^s$ diffeomorphisms are dense in the space $\mathrm{Diff}^r(M)$ and, therefore, in any Newhouse domain $\mathcal{N} \subset \mathrm{Diff}^r(M)$, see e.g. [PM].

THEOREM 8 ([GST]). *For any positive integers $s > k \geq r$ the set of $C^s$ diffeomorphisms having a k-degenerate periodic orbit is dense in a Newhouse domain $\mathcal{N} \subset \mathrm{Diff}^r(M)$.*

A proof of this theorem is outlined in [GST] (see Theorem 4). The authors told me that they are preparing a detailed proof. A slightly different detailed proof can also be found in [K1].

Theorem 8 and Newhouse's theorem imply the following important result:

COROLLARY 2 ([GST]). *Let $f \in \mathrm{Diff}^r(M)$ be a diffeomorphism exhibiting a homoclinic tangency. There is no finite number $s$ such that a generic s-parameter family $\{f_\varepsilon\}$ unfolding a diffeomorphism $f_0 = f$ is a versal family of $f_0$, meaning that the family $\{f_\varepsilon\}$ describes all possible bifurcations occurring next to $f$. Indeed, to describe all possible bifurcations of a k-degenerate periodic orbit one needs at least $k + 1$ parameters and $k$ can be arbitrarily large.*

4.1. *A Proof of Theorem 6.* Fix a $C^r$-metric $\rho_r$ in $\mathrm{Diff}^r(M)$ defined in the standard way (see e.g. [PM]). Let $f$ be a $C^r$ diffeomorphism which belongs to a Newhouse domain $\mathcal{N}$. Write $f \mapsto_{\varepsilon,r} g$ if $g$ is a $C^r$-perturbation of size at most $\varepsilon$ with respect to $\rho_r$. Consider an arbitrary sequence of positive integer numbers $\{a_n\}_{n=1}^\infty$.

Now for any $\varepsilon$ we construct a $3\varepsilon$ perturbation $f_3$ of a diffeomorphism $f$ such that for some $n_1$ the diffeomorphism $f_3$ has $n_1 a_{n_1}$ hyperbolic periodic orbits of period $n_1$. Hyperbolicity implies that the same is true for all diffeomorphisms sufficiently close to $f_4$.

Step 1. $f \mapsto_{\varepsilon,r} f_1$, where $f_1$ belongs to a Newhouse domain and is $C^\infty$ smooth.

Step 2. By Theorem 8, there exists a $C^r$-perturbation $f_1 \mapsto_{\varepsilon,r} f_2$ such that $f_2$ has a $k$-degenerate periodic orbit $q$ of an arbitrarily large period, where $k \geq r$.

Step 3. Let $n_1$ be a period of the $k$-degenerate periodic orbit $q$. It is easy to show that one can find $f_2 \mapsto_{\varepsilon,r} f_3$ such that, in a small neighborhood of $q$, $f_3$ has $n_1 a_{n_1}$ *hyperbolic periodic points* of period $n_1$.

Therfore, we show that there exists a neighborhood $U \subset \mathrm{Diff}^r(M)$ which is arbitrarily $C^r$-close to $f$ with the following property for all $g \in U$:

$$(6) \qquad \frac{\#\{x : g^n(x) = x\}}{a_n} \geq n.$$

If the diffeomorphism $f_1$ belongs to a Newhouse domain $\mathcal{N} \subset \mathrm{Diff}^r(M)$, then we can choose perturbation in steps 2-4 so small that $f_3$ belongs to the same Newhouse domain $\mathcal{N}$. It is not difficult to see from steps 1–3 that for an open dense set in $\mathcal{N}$ the condition (6) holds at least for one $n$. Iterating steps 1–3 one constructs a residual set such that for each diffeomorphism $f$ from that residual set, the condition (6) holds for infinitely many $n$'s. This completes the proof of Theorem 6 for the case $\dim M = 2$.

Note that Newhouse used similar inductive argument to prove his well-known phenomenon on infinitely many coexisting sinks [N], [PT], [R], and [TY].

4.2. *A Proof of Theorem* 6 *for any* $\dim M > 2$. We shall use the construction described in Step 1 of the proof of the Main Theorem.

Consider a compact manifold $M$ of dimension $\dim M > 2$ and a diffeomorphism $F \in \mathrm{Diff}^r(M)$. Fix a sequence of numbers $\{a_n\}_{n \in \mathbb{Z}_+}$. Suppose $F$ has a $C^r$-stable invariant two-dimensional manifold $N \subset M$ and the restriction diffeomorphism $f = F|_N : N \to N$ belongs to a Newhouse domain $\mathcal{N} \subset \mathrm{Diff}^r(N)$. $C^r$-*stability of the invariant manifold* $N$ means that any $C^r$-perturbation $\tilde{F} \in \mathrm{Diff}^r(M)$ of $F$ also has a two-dimensional invariant manifold $\tilde{N}$ which is $C^r$-close to $N$ and which induces a diffeomorphism $\tilde{f} = \tilde{F}|_N : N \to N$ which is $C^r$-close to the restriction $f = F|_N : N \to N$ (see Step 1 in Section 2 above for an exact formula for $\tilde{f}$). The Sacker theorem [Sa] gives an explicit condition when $F$ has a $C^r$-stable invariant manifold. It is important that this is an open condition in $\mathrm{Diff}^r(M)$.

We proved in the last subsection that the set of diffeomorphisms for which the condition (6) is satisfied for at least one $n \in \mathbb{Z}_+$ is open and dense in a Newhouse domain $\mathcal{N} \subset \mathrm{Diff}^r(N)$. This implies that in a neighborhood $U$ of $F$ in the space of diffeomorphisms $\mathrm{Diff}^r(M)$ there is an open and dense set $\mathcal{D}_1$ of diffeomorphisms such that each one satisfies the condition (6) for some $n \in \mathbb{Z}_+$.

Let $\mathcal{D}_{n_1}$ be an open subset of $U$ consisting of diffeomorphisms for which the condition (6) holds (for $n = n_1$). The union $\mathcal{D}^1 = \cup_{n_1} \mathcal{D}_{n_1}$ is open and dense in $U$.

For each $n_1 \in Z_+$ and a sufficiently large $n_2 > n_1$ there is an open subset $\mathcal{D}_{n_1,n_2}$ in $\mathcal{D}_{n_1}$ of diffeomorphisms satisfying the condition (6) (for both $n = n_1$ and $n = n_2$). Moreover, the union $\mathcal{D}^1_{n_1} = \cup_{n_2} \mathcal{D}_{n_1,n_2}$ is open and dense in $\mathcal{D}_{n_1}$. Therefore, the union $\mathcal{D}^2 = \cup_{n_1 < n_2} \mathcal{D}_{n_1,n_2}$ is an open and dense in $U$. Continuing, we may define $\mathcal{D}^3, \mathcal{D}^4, \dots$ Then $\cap_r \mathcal{D}^r$ is a residual set in the open set $U$. Moreover, any member of $\cap_r \mathcal{D}^r$ satisfies condition (6) for infinitely many of $n_i$'s. This completes the proof of Theorem 6. □

*Acknowledgments.* I would like to express warmest thanks to my thesis advisor John Mather. He proposed to look at the problem of growth of a number of periodic orbits and suggested the pertinent idea that a highly degenerate periodic orbit can generate a lot of periodic orbits in an open way. Several discussions of elimination theory with him were fruitful for me. J. Mather and G. Forni have made useful remarks on the text. A. Katok pointed out to me the question of Artin-Mazur. Special thanks to G. Levin who suggested to me an important idea of elimination of variables. S. Patinkin corrected the English of the paper. Let me express my sincere gratitude to all of them and J. Milnor for moral support which was important to me.

PRINCETON UNIVERSITY, PRINCETON NJ
*E-mail address*: kaloshin@math.princeton.edu